# SUMMATION TEST FOR GAP PENALTIES AND STRONG LAW OF THE LOCAL ALIGNMENT SCORE[1]

By Hock Peng Chan

*National University of Singapore*

A summation test is proposed to determine admissible types of gap penalties for logarithmic growth of the local alignment score. We also define a converging sequence of log moment generating functions that provide the constants associated with the large deviation rate and logarithmic strong law of the local alignment score and the asymptotic number of matches in the optimal local alignment.

**1. Introduction.** In protein and DNA sequence matching, two sequences of length $m$ and $n$ are aligned to determine if they have a segment each that is significantly matched. A local alignment score is assigned according to the quality of the matches in the alignment subtracted by penalties for gaps present within the alignment. The gap penalty is of the form $\Delta + \gamma(k)$ [with $\gamma(1) = 0$], for a gap of length $k$. The choice of $\Delta$, also known as the gap initialization penalty, reflects our belief in the frequency of segment insertion/deletion in the evolutionary process; while the choice of $\gamma$ reflects our belief in the distribution of the length of the segment that is inserted into or deleted from DNA or protein sequences.

The affine gap penalty function corresponds to $\gamma(k) = \delta(k-1)$ for some $\delta > 0$ and is currently the most popular in sequence alignment programs (cf. BLAST; [2]). Part of its popularity can be attributed to the recursive Smith–Waterman algorithm (cf. [17]) that allows the local alignment score to be computed in $O(mn)$ time (cf. [11]). Much research has been done to understand the asymptotic behavior of the local alignment score for affine gap penalties. In [3], it was shown that the gap penalties can be essentially divided into two types; according to whether the local alignment score grows at a logarithmic rate or linear rate. Logarithmic rate of growth is statistically

Received September 2003; revised September 2004.
[1]Supported in part by the National University of Singapore RP-115-000-031-112.
*AMS 2000 subject classifications.* Primary 60F10, 60F15; secondary 92D20.
*Key words and phrases.* Sequence alignment, large deviations, log moment generating function, change of measure.

desirable and, hence, the condition provided in the paper for determining logarithmic growth is useful in practice.

Under the Hidden Markov Model (HMM) theory (cf. [10]), the local alignment score for affine gap penalties is equivalent to the maximum likelihood score under the assumption that the length of segments inserted or deleted is geometrically distributed. This does not agree with extensive empirical studies (cf. [4, 12]) which show that a heavier tail distribution is more likely and suggests the use of long-range gap penalties that satisfy $\gamma(k) = o(k)$. Common gap penalties that have been considered are the power law $[\gamma(k) = \delta(k-1)^\alpha$ for some $0 < \alpha < 1]$ and the logarithmic $[\gamma(k) = \delta \log k]$ gap penalties. Algorithms using $O(mn)$ time for computing the local alignment score are available (cf. [14, 19] for global alignments and [15] for local alignments). However, there has so far been little understanding of the asymptotics of the local alignment score for these gap penalties and questions about the appropriateness of these scores for statistical analysis remains.

Over the past decade, there have been many advances in the use of alignment algorithms for the prediction of RNA secondary structure from primary sequences; see, for example, [5] and [18] for an overview of the underlying issues. The interaction energy of base pairings are used to determine the scores of similarity matrices, while unaligned regions are associated with loops, which require a logarithmic "loop energy" for their formation, and this supports the use of the logarithmic gap penalty function. The superior performance of the logarithmic and power law gap penalty functions in deriving biologically meaningful optimal alignments was confirmed in a detailed study by Dewey [9]. In the alignment of weakly related proteins, it was also observed that long intervening loops are relatively nonconserved and best left unaligned (cf. [1]). Long-range penalty functions, which encourages long gaps, are suitable for this purpose.

In Section 2, we provide a simple summation test for $\gamma$ that can determine if logarithmic growth of the local alignment score is possible. Section 3 extends these results into a strong law under the additional assumption that $\lim_{k \to \infty} \gamma(k)/\log k = \infty$. In Section 4 the asymptotic number of matches in the optimal alignment was also shown to obey a strong law.

**2. A summation test for gap penalties.** Let $\mathcal{A}$ be a finite alphabet which can be used to represent either the twenty amino-acids in protein sequences or the four nucleotide bases in DNA sequences. Let $K : \mathcal{A} \times \mathcal{A} \to \mathbf{R}$ be a similarity score matrix satisfying $K(a,b) = K(b,a)$ for all $a, b \in \mathcal{A}$ and let $g : \{0, 1, \ldots\} \to [0, \infty)$ with $g(0) = 0$ be a nondecreasing, concave [i.e., $g(k+1) - g(k) \leq g(k) - g(k-1)$] function. Let $\mathcal{Z}$ be the class of all nonempty candidate alignments $\mathbf{z} = \{(i(t), j(t)) : 1 \leq t \leq u\}$, where $i(1) < \cdots < i(u)$ and $j(1) < \cdots < j(u)$ are positive integers. We shall use the notation $\mathbf{z}(u)$ to signify a candidate alignment with $u$ pairs or matches. Throughout the text, $|\cdot|$ shall denote the number of elements in a finite set and, in particular, $|\mathbf{z}(u)| = u$. Given sequences $\mathbf{x}_m = x_1 \cdots x_m$ and $\mathbf{y}_n = y_1 \cdots y_n$, where $x_i, y_j \in \mathcal{A}$ for all $i$ and $j$, we define

$$S_{\mathbf{z}(u)}(\mathbf{x}_m, \mathbf{y}_n) = \sum_{t=1}^{u} K(x_{i(t)}, y_{j(t)})$$
$$- \sum_{t=1}^{u-1} [g(i(t+1) - i(t) - 1) + g(j(t+1) - j(t) - 1)] \quad (2.1)$$

if $i(u) \leq m$ and $j(u) \leq n$. For completeness, define $S_{\mathbf{z}(u)}(\mathbf{x}_m, \mathbf{y}_n) = -\infty$ if either $i(u) > m$ or $j(u) > n$. Let the local alignment score

$$H(\mathbf{x}_m, \mathbf{y}_n) = \sup_{\mathbf{z} \in \mathcal{Z}} S_{\mathbf{z}}(\mathbf{x}_m, \mathbf{y}_n). \quad (2.2)$$

Under the null hypothesis of no relation between $\mathbf{x}_m$ and $\mathbf{y}_n$, we assume that $x_1, x_2, \ldots$ and $y_1, y_2, \ldots$ are independent and identically distributed with probability measure $\mu$ satisfying $\mu(a) > 0$ for all $a \in \mathcal{A}$. Define $\mu(\mathbf{x}_m) = \prod_{i=1}^{m} \mu(x_i)$ and assume that

$$E[K(x_1, y_1)] < 0, \qquad K_{\max} := \max_{a,b \in \mathcal{A}} K(a,b) > 0. \quad (2.3)$$

The local alignment score for gapless alignments, denoted by $H_\infty$, can be expressed in the form (2.1)–(2.2) by setting $g(k) = \infty$ for all $k \geq 1$. Its asymptotic behavior was studied in [7, 8]. Let $\theta_*$ be the unique positive solution to the equation $E \exp[\theta K(x_1, y_1)] = 1$. It was shown that under (2.3), $H_\infty(\mathbf{x}_n, \mathbf{y}_n)$ has an asymptotic Gumbel distribution and

$$H_\infty(\mathbf{x}_n, \mathbf{y}_n) / \log n \to 2/\theta_* \quad \text{a.s. as } n \to \infty. \quad (2.4)$$

Analogous to (2.1)–(2.2), we define

$$R_{\mathbf{z}(u)}(\mathbf{x}_m, \mathbf{y}_n) = S_{\mathbf{z}(u)}(\mathbf{x}_m, \mathbf{y}_n) - g(i(1) - 1)$$
$$- g(j(1) - 1) - g(m - i(u)) - g(n - j(u)), \quad (2.5)$$

$$G(\mathbf{x}_m, \mathbf{y}_n) = \sup_{\mathbf{z} \in \mathcal{Z}} R_{\mathbf{z}}(\mathbf{x}_m, \mathbf{y}_n). \quad (2.6)$$

$G$ is known as the global alignment score and differs from the local alignment score $H$ in that unaligned letters both before and after the alignment $\mathbf{z}$ are penalized. If $g$ is chosen such that $\beta := \lim_{n \to \infty} E[G(\mathbf{x}_n, \mathbf{y}_n)/n] > 0$, then $H(\mathbf{x}_n, \mathbf{y}_n)/n \to \beta$ in probability and the gap penalty is said to lie in the linear domain. Conversely, for $\beta < 0$, there exists $\tau_2 > \tau_1 > 0$ such that $\lim_{n \to \infty} P\{\tau_2 > H(\mathbf{x}_n, \mathbf{y}_n)/\log n > \tau_1\} = 1$ and the gap penalty is said to lie in the logarithmic domain (cf. [3], Lemmas 2 and 3).

In some sequence alignment software, for example, XPARAL (cf. [13]), the user is required to specify gap penalties of the form $g(k) = \Delta + \gamma(k)$ for $k \geq 1$, where $\gamma(k) = \delta \log k$, $\gamma(k) = \delta(k-1)^{1/2}$ or $\gamma(k) = \delta(k-1)$ for some $\delta > 0$. By Arratia and Waterman [3], it follows that if $\gamma(k) = \delta(k-1)$, then $g$ lies in the logarithmic region if the gap penalty is chosen large enough. However, it is unclear for the cases $\gamma(k) = \delta \log k$ and $\gamma(k) = \delta(k-1)^{1/2}$ that logarithmic growth of $H$ is possible. Note that for these choices of $\gamma$, the constant $\beta$ is always nonnegative. This can be seen by considering an alignment with exactly one match. In Theorem 1, we provide a summation test that will allow us to determine the types of $\gamma$ for which logarithmic growth occurs when $\Delta$ is large. It formalizes a statement in [16], where a rough heuristic is used to suggest that gap penalties satisfying $\sum_{k=1}^{\infty} \exp[-\theta_* g(k)] < \infty$ should be chosen for logarithmic growth of $H$.

THEOREM 1. *Let $g(k) = \Delta + \gamma(k)$ for $k \geq 1$ with $\gamma(1) = 0$.*

(a) *If $\sum_{k=1}^{\infty} \exp[-\widehat{\theta}\gamma(k)] < \infty$ for some $\widehat{\theta} < \theta_*$, then $g$ lies in the logarithmic domain for all large $\Delta$.*

(b) *If $\sum_{k=1}^{\infty} \exp[-\widehat{\theta}\gamma(k)] = \infty$ for some $\widehat{\theta} > \theta_*$, then $g$ lies in the linear domain for all $\Delta$.*

(c) *Let $\gamma(k) = \delta \log k$ for some $\delta > 0$. If $\delta > \theta_*^{-1}$, then $g$ lies in the logarithmic domain for all large $\Delta$. Conversely, if $\delta < \theta_*^{-1}$, then $g$ lies in the linear domain for all $\Delta$.*

Let $\mathcal{Z}_\kappa = \{\mathbf{z} \in \mathcal{Z} : (1,1) \in \mathbf{z}, |\mathbf{z}| = \kappa\}$. Define

$$(2.7) \qquad G_\kappa(\mathbf{x}_m, \mathbf{y}_n) = \max_{\zeta \in \mathcal{Z}_\kappa} R_\zeta(\mathbf{x}_m, \mathbf{y}_n).$$

For $\theta > 0$, define

$$(2.8) \qquad h_\kappa(\theta) = \sum_{m,n \geq \kappa} E \exp[\theta G_\kappa(\mathbf{x}_m, \mathbf{y}_n)], \qquad \psi_\kappa(\theta) = \log h_\kappa(\theta).$$

To prove Theorem 1, we need to consider only $\kappa = 1$, but the strong law results of Theorems 2 and 3 use the convergence of $\psi_\kappa(\theta)/\kappa$ as $\kappa \to \infty$. We preface the proof of Theorem 1 with Lemma 1, which uses an importance sampling scheme to achieve a change of measure. For $\kappa = 1$ and $g(k) = \Delta + \delta(k-1)$, a modified version of this scheme was implemented in [6] for efficient simulation of $P\{H(\mathbf{x}_m, \mathbf{y}_n) \geq c\}$.

LEMMA 1. *Let $\theta_\kappa$ be a positive root of $\psi_\kappa(\theta) = 0$ (if it exists). Then*

$$(2.9) \qquad P\{H(\mathbf{x}_m, \mathbf{y}_n) \geq c\} \leq mn \exp[\theta_\kappa(\kappa - 1)K_{\max}] \exp(-\theta_\kappa c).$$

PROOF. Let us simulate $(\mathbf{x}_m, \mathbf{y}_n)$ in the following manner:

1. Initialization step. Simulate $(i_*, j_*)$ uniformly from $\{1, \ldots, m\} \times \{1, \ldots, n\}$ and let $x_i, y_j \sim \mu$ for all $i < i_*$ and $j < j_*$. Initialize partial sum $S = 0$.

2. Repetition steps.

(a) Simulation. Simulate $(\mathbf{v}_r, \mathbf{w}_s)$ from the measure $\nu$ on the domain $\bigcup_{r=\kappa}^{\infty} \mathcal{A}^r \times \bigcup_{s=\kappa}^{\infty} \mathcal{A}^s$ with

(2.10) $$\nu(\mathbf{v}_r, \mathbf{w}_s) = \exp[\theta_\kappa G_\kappa(\mathbf{v}_r, \mathbf{w}_s)]\mu(\mathbf{v}_r)\mu(\mathbf{w}_s).$$

Note that both $r$ and $s$ are random here, taking values in $\{\kappa, \kappa+1, \dots\}$. Moreover, $\nu$ is a probability measure because $\psi_\kappa(\theta_\kappa) = 0$.

(b) Check that segment generated is not too long. If $i_* + r - 1 > m$ or $j_* + s - 1 > n$, go to step 3. Otherwise, proceed to (c).

(c) Updating. Let $x_{i_*} \cdots x_{i_*+r-1} = \mathbf{v}_r$, $y_{j_*} \cdots y_{j_*+s-1} = \mathbf{w}_s$ and let $(\text{new})(i_*, j_*) = (\text{old})(i_*, j_*) + (r, s)$. Let $(\text{new})S = (\text{old})S + G_\kappa(\mathbf{v}_r, \mathbf{w}_s)$. If $S \geq c - (\kappa-1)K_{\max}$, go to step 3. Otherwise, repeat step 2.

3. Conclusion step. Simulate $x_i, y_j \sim \mu$ for all $i \geq i_*$ and $j \geq j_*$.

Let $Q$ denote the probability measure of $(\mathbf{x}_m, \mathbf{y}_n)$ simulated in this manner and let $P(\mathbf{x}_m, \mathbf{y}_n) = \mu(\mathbf{x}_m)\mu(\mathbf{y}_n)$. Equation (2.9) clearly holds when $c \leq (\kappa-1)K_{\max}$ so we may assume without loss of generality that $c > (\kappa-1)K_{\max}$. Let $A = \{(\mathbf{x}_m, \mathbf{y}_n) : H(\mathbf{x}_m, \mathbf{y}_n) \geq c\}$. For all $(\mathbf{x}_m, \mathbf{y}_n) \in A$, there exists $\mathbf{z} \in \mathcal{Z}$ such that $S_\mathbf{z}(\mathbf{x}_m, \mathbf{y}_n) \geq c$. Since $c > (\kappa-1)K_{\max}$, it follows that $\mathbf{z}$ has at least $\kappa$ matches and can be expressed in the form $\mathbf{z} = \{(i(t), j(t)) : 1 \leq t \leq \lambda\kappa + q\}$ for some $\lambda \geq 1$ and $0 \leq q < \kappa$. Let $\zeta = \{(i(t), j(t)) : 1 \leq t \leq \lambda\kappa\}$, which is $\mathbf{z}$ without the last $q$ matches. Since $q \leq (\kappa-1)$ and $S_\mathbf{z}(\mathbf{x}_m, \mathbf{y}_n) \geq c$, it follows that

(2.11) $$S_\zeta(\mathbf{x}_m, \mathbf{y}_n) \geq c - (\kappa-1)K_{\max}.$$

We break-up $x_{i(1)}x_{i(1)+1} \cdots x_{i(\lambda\kappa)}$ into $\lambda$ segments $\mathbf{v}^{(1)}, \dots, \mathbf{v}^{(\lambda)}$ with $\mathbf{v}^{(1)} = x_{i(1)} \cdots x_{i(\kappa+1)-1}$, $\mathbf{v}^{(2)} = x_{i(\kappa+1)} \cdots x_{i(2\kappa+1)-1}, \dots, \mathbf{v}^{(\lambda-1)} = x_{i((\lambda-2)\kappa+1)} \cdots x_{i((\lambda-1)\kappa+1)-1}$ and for the last segment, $\mathbf{v}^{(\lambda)} = x_{i((\lambda-1)\kappa+1)} \cdots x_{i(\lambda\kappa)}$. Similarly, $y_{j(1)}y_{j(1)+1} \cdots y_{j(\lambda\kappa)}$ is broken up into $\lambda$ segments $\mathbf{w}^{(1)}, \dots, \mathbf{w}^{(\lambda)}$ where $\mathbf{w}^{(1)} = y_{j(1)} \cdots y_{j(\kappa+1)-1}$, $\mathbf{w}^{(2)} = y_{j(\kappa+1)} \cdots y_{j(2\kappa+1)-1}, \dots$ and for the last segment $\mathbf{w}^{(\lambda)} = y_{j((\lambda-1)\kappa+1)} \cdots y_{j(\lambda\kappa)}$. $(\mathbf{x}_m, \mathbf{y}_n)$ can be generated from the simulation scheme above if $(i(1), j(1))$ is simulated in step 1 [as $(i_*, j_*)$] and $(\mathbf{v}^{(\eta)}, \mathbf{w}^{(\eta)})$ are generated on the $\eta$th iteration of step 2(a). Since $S_\zeta(\mathbf{x}_m, \mathbf{y}_n) \leq \sum_{\eta=1}^{\lambda} G_\kappa(\mathbf{v}^{(\eta)}, \mathbf{w}^{(\eta)})$, it follows by (2.10) and (2.11) that

(2.12)
$$\frac{dQ}{dP}(\mathbf{x}_m, \mathbf{y}_n) \geq (nm)^{-1} \prod_{\eta=1}^{\lambda} [\nu(\mathbf{v}^{(\eta)}, \mathbf{w}^{(\eta)})/\mu(\mathbf{v}^{(\eta)})\mu(\mathbf{w}^{(\eta)})]$$
$$\geq (nm)^{-1} \exp[\theta_\kappa S_\zeta(\mathbf{x}_m, \mathbf{y}_n)]$$
$$\geq (nm)^{-1} \exp\{\theta_\kappa[c - (\kappa-1)K_{\max}]\}.$$

This holds for all $(\mathbf{x}_m, \mathbf{y}_n) \in A$ and, hence, (2.9) follows from the identity $P(A) = E[(dP/dQ)\mathbf{1}_A]$. □

PROOF OF THEOREM 1. Since $\mathcal{Z}_1$ contains only the alignment $\{(1,1)\}$, it follows that $G_1(\mathbf{x}_m, \mathbf{y}_n) = K(x_1, y_1) - g(m-1) - g(n-1)$ and, hence, by (2.8),

$$(2.13) \qquad h_1(\theta) = \left[1 + e^{-\theta \Delta} \sum_{k=1}^{\infty} e^{-\theta \gamma(k)}\right]^2 E \exp[\theta K(x_1, y_1)].$$

Let $\sum_{k=1}^{\infty} \exp[-\widehat{\theta}\gamma(k)] < \infty$ for some $0 < \widehat{\theta} < \theta_*$. Since $E \exp[\widehat{\theta} K(x_1, y_1)] < 1$, we can find $\Delta$ large enough such that $h_1(\widehat{\theta}) < 1$. By (2.3) and (2.13), $\psi_1(\theta) = \log h_1(\theta) \to \infty$ as $\theta \to \infty$ and, hence, $\psi_1(\theta_1) = 0$ for some $\theta_1 > \widehat{\theta}$. By Lemma 1 with $\kappa = 1$,

$$(2.14) \qquad \lim_{n \to \infty} P\{H(\mathbf{x}_n, \mathbf{y}_n) \geq 3\log n/\theta_1\} = 0.$$

Since $H \geq H_\infty$, the gapless local alignment score, it follows from (2.4) that $\lim_{n \to \infty} P\{H(\mathbf{x}_n, \mathbf{y}_n) \geq \log n/\theta_*\} = 1$ and, hence, (a) follows from (2.14). The first part of (c) also follows from (a) by choosing $\widehat{\theta} \in (\delta^{-1}, \theta_*)$.

We shall next show the second part of (c). Let $g(k) = \Delta + \delta \log k$ for some $\delta < \theta_*^{-1}$. Let $\mathbf{v}^{(\eta)} = x_{r(\eta-1)+1} \cdots x_{r\eta}$ and $\mathbf{w}^{(\eta)} = y_{r(\eta-1)+1} \cdots y_{r\eta}$ for $1 \leq \eta \leq \lambda$, where $\lambda$ and $r$ are positive integers to be specified later. Then $G(\mathbf{x}_{r\lambda}, \mathbf{y}_{r\lambda}) \geq \sum_{\eta=1}^{\lambda} H_\infty(\mathbf{v}^{(\eta)}, \mathbf{w}^{(\eta)}) - 2(\lambda+1)g(2r)$ and, hence, it follows from (2.4) that for any $\varepsilon > 0$, there exists $r$ large enough such that

$$(2.15) \begin{aligned} E[G(\mathbf{x}_{r\lambda}, \mathbf{y}_{r\lambda})] &\geq \lambda E[H_\infty(\mathbf{x}_r, \mathbf{y}_r)] - 2(\lambda+1)g(2r) \\ &\geq 2\lambda(1-\varepsilon)(\log r)/\theta_* - 2(\lambda+1)[\Delta + \delta \log(2r+1)]. \end{aligned}$$

Since $\delta < \theta_*^{-1}$, it follows from (2.15) that there exists $\varepsilon$ small enough and $\lambda, r$ large enough such that $E[G(\mathbf{x}_{r\lambda}, \mathbf{y}_{r\lambda})] > 0$ and, hence, $g$ lies in the linear domain.

To show (b), pick $\delta \in (\widehat{\theta}^{-1}, \theta_*^{-1})$. Since $\sum_{k=1}^{\infty} \exp[-\widehat{\theta}(\delta \log k)] < \infty$ and $\sum_{k=1}^{\infty} \exp[-\widehat{\theta}\gamma(k)] = \infty$, it follows that $\gamma(k) \leq \delta \log k$ for infinitely many $k$. Then for any $\varepsilon > 0$, (2.15) holds for infinitely many $r$ and (b) follows by choosing $\lambda, r$ large enough and $\varepsilon$ small enough. □

**3. Large deviations and the strong law of large numbers.** In this section we extend the large deviations and strong law results of Arratia and Waterman [3] and Zhang [20] to gap penalties satisfying $\lim_{k \to \infty} g(k)/k = 0$, by-passing the Azuma–Hoeffding inequality that was central to their proofs for the case $\lim_{k \to \infty} g(k)/k > 0$.

THEOREM 2. *Let $g(k) = \Delta + \gamma(k)$ for $k \geq 1$, where $\gamma(1) = 0$ and $\lim_{k \to \infty} \gamma(k)/\log k = \infty$. Then $\psi(\theta) = \lim_{\kappa \to \infty} \psi_\kappa(\theta)/\kappa$ is well defined, convex and finite for all $\theta > 0$. Moreover, for all large $\Delta$,*

$$(3.1) \qquad \psi(\theta) = 0 \text{ has a unique positive root } \widetilde{\theta}.$$

*Under (3.1), the following also holds.*

(a) *If* $\min\{m,n\}/c \to \infty$ *and* $\log(mn) = o(c)$ *as* $c \to \infty$, *then*

(3.2) $$\lim_{c \to \infty} -c^{-1} \log P\{H(\mathbf{x}_m, \mathbf{y}_n) \geq c\} = \widetilde{\theta}.$$

(b) $H(\mathbf{x}_n, \mathbf{y}_n)/\log n \to 2/\widetilde{\theta}$ *a.s. as* $n \to \infty$.

PROOF. Let $\kappa, \eta$ be positive integers and consider $\mathbf{x}_m, \mathbf{y}_n$ with $m, n \geq \kappa + \eta$. Let

(3.3) $$\Pi(\mathbf{x}_m, \mathbf{y}_n) = \{(\mathbf{v}_p^{(1)}, \mathbf{v}_q^{(2)}, \mathbf{w}_r^{(1)}, \mathbf{w}_s^{(2)}) : \mathbf{x}_m = \mathbf{v}_p^{(1)} \mathbf{v}_q^{(2)}, \mathbf{y}_n = \mathbf{w}_r^{(1)} \mathbf{w}_s^{(2)}$$
$$\text{with } p, r \geq \kappa \text{ and } q, s \geq \eta\}.$$

In other words, $\mathbf{v}_p^{(1)} = x_1 \cdots x_p$, $\mathbf{v}_q^{(2)} = x_{p+1} \cdots x_m$, $\mathbf{w}_r^{(1)} = y_1 \cdots y_r$ and $\mathbf{w}_s^{(2)} = y_{r+1} \cdots y_n$. For notational simplicity, we shall henceforth omit superscripts (1), (2) when describing members of $\Pi(\mathbf{x}_m, \mathbf{y}_n)$. Since $G_{\kappa+\eta}(\mathbf{x}_m, \mathbf{y}_n) = S_{\mathbf{z}}(\mathbf{x}_m, \mathbf{y}_n)$ for some $\mathbf{z} = \{(i(t), j(t)) : 1 \leq t \leq \kappa + \eta\} \in \mathcal{Z}_{\kappa+\eta}$, it follows by selecting $p = i(\kappa+1) - 1$ and $q = j(\kappa+1) - 1$ that $G_{\kappa+\eta}(\mathbf{x}_m, \mathbf{y}_n) = G_\kappa(\mathbf{v}_p, \mathbf{w}_r) + G_\eta(\mathbf{v}_q, \mathbf{w}_s)$ for some $(\mathbf{v}_p, \mathbf{v}_q, \mathbf{w}_r, \mathbf{w}_s) \in \Pi(\mathbf{x}_m, \mathbf{y}_n)$. As $\mu(\mathbf{x}_m) = \mu(\mathbf{v}_p)\mu(\mathbf{v}_q)$ and $\mu(\mathbf{y}_n) = \mu(\mathbf{w}_r)\mu(\mathbf{w}_s)$, the inequality

(3.4) $$\exp[\theta G_{\kappa+\eta}(\mathbf{x}_m, \mathbf{y}_n)]\mu(\mathbf{x}_m)\mu(\mathbf{y}_n)$$
$$\leq \sum_{(\mathbf{v}_p, \mathbf{v}_q, \mathbf{w}_r, \mathbf{w}_s) \in \Pi(\mathbf{x}_m, \mathbf{y}_n)} \exp[\theta G_\kappa(\mathbf{v}_p, \mathbf{w}_r)]\mu(\mathbf{v}_p)\mu(\mathbf{w}_r)$$
$$\times \exp[\theta G_\eta(\mathbf{v}_q, \mathbf{w}_s)]\mu(\mathbf{v}_q)\mu(\mathbf{w}_s)$$

holds because there exists a term on the right-hand side of (3.4) that is equal to the left-hand side. Noting that

(3.5) $$\bigcup_{(\mathbf{x}_m, \mathbf{y}_n) : m, n \geq \kappa + \eta} \Pi(\mathbf{x}_m, \mathbf{y}_n) = \{(\mathbf{v}_p, \mathbf{v}_q, \mathbf{w}_r, \mathbf{w}_s) : p, r \geq \kappa \text{ and } q, s \geq \eta\}$$

and that the left-hand side of (3.5) is a disjoint union, we can conclude from (3.4) that, for all $\theta > 0$,

$$\psi_{\kappa+\eta}(\theta) = \log\left\{\sum_{(\mathbf{x}_m, \mathbf{y}_n) : m, n \geq \kappa + \eta} \exp[\theta G_{\kappa+\eta}(\mathbf{x}_m, \mathbf{y}_n)]\mu(\mathbf{x}_m)\mu(\mathbf{y}_n)\right\}$$

$$\leq \log\left\{\sum_{(\mathbf{x}_m, \mathbf{y}_n) : m, n \geq \kappa + \eta} \left[\sum_{(\mathbf{v}_p, \mathbf{v}_q, \mathbf{w}_r, \mathbf{w}_s) \in \Pi(\mathbf{x}_m, \mathbf{y}_n)} \exp[\theta G_{\kappa+\eta}(\mathbf{v}_p, \mathbf{w}_r)]\right.\right.$$
$$\times \mu(\mathbf{v}_p)\mu(\mathbf{w}_r)$$
(3.6) $$\left.\left.\times \exp[\theta G_\eta(\mathbf{v}_q, \mathbf{w}_s)]\mu(\mathbf{v}_q)\mu(\mathbf{w}_s)\right]\right\}$$

$$= \log \Biggl\{ \sum_{(\mathbf{v}_p, \mathbf{v}_q, \mathbf{w}_r, \mathbf{w}_s):\, p, r \geq \kappa \text{ and } q, s \geq \eta} \exp[\theta G_\kappa(\mathbf{v}_p, \mathbf{w}_r)] \mu(\mathbf{v}_p) \mu(\mathbf{w}_r)$$
$$\times \exp[\theta G_\eta(\mathbf{v}_q, \mathbf{w}_s)] \mu(\mathbf{v}_q) \mu(\mathbf{w}_s) \Biggr\}$$

$$= \log \Biggl\{ \sum_{(\mathbf{v}_p, \mathbf{w}_r):\, p \geq \kappa \text{ and } r \geq \eta} \exp[\theta G_\kappa(\mathbf{v}_p, \mathbf{w}_r)] \mu(\mathbf{v}_p) \mu(\mathbf{w}_r)$$
$$\times \sum_{(\mathbf{v}_q, \mathbf{w}_s):\, q \geq \kappa \text{ and } s \geq \eta} \exp[\theta G_\eta(\mathbf{v}_q, \mathbf{w}_s)] \mu(\mathbf{v}_q) \mu(\mathbf{w}_s) \Biggr\}$$

$$= \psi_\kappa(\theta) + \psi_\eta(\theta).$$

Moreover, as $h_\kappa(\theta) \geq E \exp[\theta G_\kappa(\mathbf{x}_\kappa, \mathbf{y}_\kappa)] = \{E \exp[\theta K(x_1, y_1)]\}^\kappa$ and $\psi_\kappa(\theta) = \log h_\kappa(\theta)$, it follows that

(3.7) $\quad \psi_\kappa(\theta)/\kappa \geq \log E \exp[\theta K(x_1, y_1)] > -\infty \quad$ for all $\theta > 0$ and $\kappa \geq 1$.

The subadditive property (3.6) then ensures that $\psi(\theta) = \lim_{\kappa \to \infty} \psi_\kappa(\theta)/\kappa$ is well defined and finite. Since $\psi_\kappa(\theta)$ is convex for all $\kappa$ (see Section A.1), it follows that $\psi$ is convex and continuous.

Pick an arbitrary positive $\widehat{\theta} < \theta_*$, where $\theta_*$ is the unique positive root of the equation $E \exp[\theta K(x, y)] = 1$. By (2.13), $\psi_1(\widehat{\theta}) < 0$ for all large $\Delta$ and, hence, $\psi(\widehat{\theta}) \leq \psi_1(\widehat{\theta}) < 0$. By (2.3) and (3.7), $\lim_{\theta \to \infty} \psi(\theta) = \infty$ and, hence, a positive solution $\widetilde{\theta}$ ($> \widehat{\theta}$) of the equation $\psi(\theta) = 0$ exists. To show that $\widetilde{\theta}$ is unique, it suffices from the convexity of $\psi$ to show that $\lim_{\theta \to 0} \psi(\theta) = 0$. Since

(3.8) $\quad G_\kappa(\mathbf{x}_r, \mathbf{y}_s) \leq \kappa K_{\max} - g(r - \kappa) - g(s - \kappa),$

it follows that

$$h_\kappa(\theta) \leq \sum_{r, s \geq \kappa} \exp\{\theta[\kappa K_{\max} - g(r - \kappa) - g(s - \kappa)]\}$$
$$= \exp(\theta \kappa K_{\max}) \Biggl\{ \sum_{k \geq 0} \exp[-\theta g(k)] \Biggr\}^2$$

and, indeed, $\psi(\theta) = \lim_{\kappa \to \infty} [\log h_\kappa(\theta)]/\kappa \leq \theta K_{\max} \to 0$ as $\theta \to 0$. Moreover, by (3.7), $\psi(\theta) = \lim_{\kappa \to \infty} \psi_\kappa(\theta)/\kappa \geq \log E \exp[\theta K(x_1, y_1)] \to 0$ as $\theta \to 0$.

(a) It follows from (3.1) that there exists $\theta_\kappa \to \widetilde{\theta}$ such that $\psi_\kappa(\theta_\kappa) = 0$ for all large $\kappa$. By Lemma 1 and as $c^{-1} \log(mn) \to 0$, it follows that

(3.9) $\quad \liminf_{c \to \infty} -c^{-1} \log P\{H(\mathbf{x}_m, \mathbf{y}_n) \geq c\} \geq \theta_\kappa.$

To get the opposite inequality, define

$$\xi_\kappa(\theta) = \sup_{r \geq \kappa} \log\{E \exp[\theta G_\kappa(\mathbf{x}_r, \mathbf{y}_r)]\}. \tag{3.10}$$

Clearly, $G_{\kappa+\eta}(\mathbf{x}_{r+s}, \mathbf{y}_{r+s}) \geq G_\kappa(\mathbf{x}_r, \mathbf{y}_s) + G_\eta(x_{r+1} \cdots x_{r+s}, y_{r+1} \cdots y_{r+s})$ and, hence,

$$E \exp[\theta G_{\kappa+\eta}(\mathbf{x}_{r+s}, \mathbf{y}_{r+s})] \geq E \exp[\theta G_\kappa(\mathbf{x}_r, \mathbf{y}_r)] E \exp[\theta G_\eta(\mathbf{x}_s, \mathbf{y}_s)].$$

By taking supremum over $r$ and $s$, the superadditive property $\xi_{\kappa+\eta}(\theta) \geq \xi_\kappa(\theta) + \xi_\eta(\theta)$ holds. Since $\psi_\kappa(\theta) \geq \xi_\kappa(\theta)$ for all $\kappa$, it follows that $\psi(\theta) \geq \lim_{\kappa \to \infty} \xi_\kappa(\theta)/\kappa$. It shall be shown in Section A.2 that if $g(k)/\log k \to \infty$, then

$$\psi(\theta) = \lim_{\kappa \to \infty} \xi_\kappa(\theta)/\kappa \qquad \text{for all } \theta > 0 \tag{3.11}$$

and, hence, there exists $\widehat{\theta}_\kappa \to \widetilde{\theta}$ such that

$$\xi_\kappa(\widehat{\theta}_\kappa) = 0 \tag{3.12}$$

for all large $\kappa$. Let $\kappa$ satisfy (3.12). It follows from (3.8) that $E \exp[\theta G_\kappa(\mathbf{x}_r, \mathbf{y}_r)] \to 0$ as $r \to \infty$, and, hence, for all $\theta > 0$, the supremum in (3.10) is attained at some $r \geq \kappa$. By (3.12), it follows that

$$E \exp[\widehat{\theta}_\kappa G_\kappa(\mathbf{x}_r, \mathbf{y}_r)] = 1 \qquad \text{for some } r \ (= r_\kappa). \tag{3.13}$$

Let $v_\eta = G_\kappa(x_{(r-1)\eta+1} \cdots x_{r\eta}, y_{(r-1)\eta+1} \cdots y_{r\eta})$ and let $Q$ be the measure under which $v_1, v_2, \ldots$ are independent with $Q\{v_\eta = k\} = \exp(\widehat{\theta}_\kappa k) P\{v_\eta = k\}$. By (3.13), $Q$ is a probability measure. Let $T_c = \inf\{\ell : \sum_{\eta=1}^\ell v_\eta \geq c\}$. Then $(dQ/dP)(v_1, \ldots, v_{T_c}) = \exp(\widehat{\theta}_\kappa \sum_{\eta=1}^{T_c} v_\eta) \leq \exp[\widehat{\theta}_\kappa (c + \kappa K_{\max})]$ whenever $T_c < \infty$. Hence,

$$P\{T_c \leq \lambda\} \geq \exp[-\widehat{\theta}_\kappa(c + \kappa K_{\max})] Q\{T_c \leq \lambda\} \tag{3.14}$$

for any positive integer $\lambda$. Pick $\lambda = \lfloor \min\{m, n\}/r \rfloor$, where $\lfloor \cdot \rfloor$ denotes the greatest integer function. Since $\min\{m, n\}/c \to \infty$, it follows that $\lambda/c \to \infty$. By the law of large numbers and as $E_Q v_1 > 0$, we can conclude that $Q\{T_c \leq \lambda\} \to 1$. By (3.14) and as $H(\mathbf{x}_m, \mathbf{y}_n) \geq \sum_{\eta=1}^\lambda v_\eta$ so that $\{H(\mathbf{x}_m, \mathbf{y}_n) \geq c\} \supset \{T_c \leq \lambda\}$, it folllows that

$$\limsup_{c \to \infty} -c^{-1} \log P\{H(\mathbf{x}_m, \mathbf{y}_n) \geq c\} \leq \widehat{\theta}_\kappa. \tag{3.15}$$

(a) then follows from (3.9) and (3.15) by letting $\kappa \to \infty$.

(b) Let $\varepsilon > 0$ and select $\kappa$ large enough such that $\psi_\kappa(\theta) = 0$ has a positive solution $\theta_\kappa$. Select a subsequence $n_k = \lfloor k^{2/\varepsilon} \rfloor + 1$. Then by Lemma 1,

$$\sum_{k=1}^\infty P\{H(\mathbf{x}_{n_k}, \mathbf{y}_{n_k}) \geq (2 + \varepsilon)(\log n_k)/\theta_\kappa\} \leq \exp[\theta_\kappa(\kappa - 1)K_{\max}] \sum_{k=1}^\infty n_k^{-\varepsilon} < \infty.$$

By the Borel–Cantelli lemma, it follows that $\limsup_{k\to\infty} H(\mathbf{x}_{n_k},\mathbf{y}_{n_k})/\log n_k \leq (2+\varepsilon)/\theta_\kappa$ a.s. Since $\log n_{k+1}/\log n_k \to 1$ and $H(\mathbf{x}_n,\mathbf{y}_n)$ is nondecreasing in $n$, it follows by choosing $\varepsilon$ arbitrarily small that

$$(3.16) \qquad \limsup_{n\to\infty} H(\mathbf{x}_n,\mathbf{y}_n)/\log n \leq 2/\theta_\kappa \qquad \text{a.s.}$$

Let $\kappa$ satisfy (3.12) and define a score matrix $\widetilde{K}$ on $\mathcal{B} := \mathcal{A}^r \ (= \mathcal{A}^{r_\kappa})$ [see (3.13)] by setting $\widetilde{K}(\mathbf{x}_r,\mathbf{y}_r) = G_\kappa(\mathbf{x}_r,\mathbf{y}_r)$. Let $\mathbf{x}^{(\eta)} = x_{(\eta-1)r+1}\cdots x_{\eta r}$ and $\mathbf{y}^{(\eta)} = y_{(\eta-1)r+1}\cdots y_{\eta r}$ for all $1 \leq \eta \leq \lambda := \lfloor n/r \rfloor$. Let $\widetilde{H}_\infty(\mathbf{x}_{r\lambda},\mathbf{y}_{r\lambda}) = \widetilde{H}_\infty(\mathbf{x}^{(1)}\cdots\mathbf{x}^{(\lambda)},\mathbf{y}^{(1)}\cdots\mathbf{y}^{(\lambda)})$ be the gapless local alignment score which treats $\mathbf{x}^{(\eta)}, \mathbf{y}^{(\eta)}$ as letters of $\mathcal{B}$ and uses $\widetilde{K}$ as the score matrix. By (3.13), it follows that (2.4) holds with $\widetilde{H}_\infty$ in place of $H_\infty$ and $\widehat{\theta}_\kappa$ in place of $\theta_*$. Hence,

$$(3.17) \quad \liminf_{n\to\infty} H(\mathbf{x}_n,\mathbf{y}_n)/\log n \geq \lim_{\lambda\to\infty} \widetilde{H}_\infty(\mathbf{x}_{r\lambda},\mathbf{y}_{r\lambda})/\log\lambda = 2/\widehat{\theta}_\kappa \qquad \text{a.s.}$$

(b) follows from (3.16) and (3.17) by letting $\kappa \to \infty$. □

**4. Asymptotic number of matches in the optimal local alignment.** For given sequences $\mathbf{x}_n, \mathbf{y}_n$, let $\mathbf{z}$ be a candidate alignment satisfying

$$(4.1) \qquad S_\mathbf{z}(\mathbf{x}_n,\mathbf{y}_n) = H(\mathbf{x}_n,\mathbf{y}_n).$$

The alignment $\mathbf{z}$ is not unique in general, but to be specific, we shall assume that there exists an ordering of the candidate alignments in $\mathcal{Z}$ and only the smallest alignment $\mathbf{z}$ with respect to this ordering that satisfies (4.1) shall be designated as the optimal local alignment and denoted by $\mathbf{z}_*$. Properties of the optimal local alignment are less stable than the local alignment score because a slight perturbation of the sequences, for example, changing one of the letters $x_i$ or $y_j$, can result in a very different optimal local alignment.

In this section our objective is to study $|\mathbf{z}_*|$, the number of matches in the optimal alignment $\mathbf{z}_*$. We shall show in Theorem 3 that under the assumptions of Theorem 2, $|\mathbf{z}_*| \sim 2\log n/\widetilde{\theta}\psi'(\widetilde{\theta})$ as $n \to \infty$ whenever the derivative $\psi'(\widetilde{\theta})$ exists. Since $H(\mathbf{x}_n,\mathbf{y}_n) \sim 2\log n/\widetilde{\theta}$ by Theorem 2(b), this gives rise to the interpretation of $\psi'(\widetilde{\theta})$ as the asymptotic score per match of the optimal alignment. The convexity of $\psi$ ensures that $\psi'(\theta)$ exists with the exception of countably many $\theta$. A more detailed discussion of the existence of $\psi'(\widetilde{\theta})$, involving measure theoretic issues, is dealt with in Section A.3.

THEOREM 3. *Let* $\lim_{k\to\infty} g(k)/\log k = \infty$ *and assume* (3.1) *holds. If* $\psi'(\widetilde{\theta})$ *is well defined, then* $|\mathbf{z}_*|/\log n \to 2/\widetilde{\theta}\psi'(\widetilde{\theta})$ *a.s.*

PROOF. Let $K_\lambda$ be a score matrix satisfying $K_\lambda(a,b) = K(a,b) + \lambda$ for all $a,b \in \mathcal{A}$. A superscript $\lambda$ in any notation defined previously will now be used to signify that the score matrix $K_\lambda$ is used. If no superscript is used, it is understood that $\lambda = 0$. Since $G_\kappa^{(\lambda)}(\mathbf{x}_m, \mathbf{y}_n) = G_\kappa(\mathbf{x}_m, \mathbf{y}_n) + \lambda\kappa$ for all $(\mathbf{x}_m, \mathbf{y}_n)$, it follows from (2.8) that $\psi_\kappa^{(\lambda)}(\theta) = \psi_\kappa(\theta) + \lambda\kappa\theta$ and, hence,

$$\psi^{(\lambda)}(\theta) = \lim_{\kappa \to \infty} \psi_\kappa^{(\lambda)}(\theta)/\kappa = \psi(\theta) + \lambda\theta. \tag{4.2}$$

By (4.2), $\psi(\widetilde{\theta}^{(\lambda)}) + \lambda\widetilde{\theta}^{(\lambda)} = \psi^{(\lambda)}(\widetilde{\theta}^{(\lambda)}) = 0 = \psi(\widetilde{\theta})$. Since $\widetilde{\theta}^{(\lambda)} \to \widetilde{\theta}$ as $\lambda \to 0$, it follows that $\psi'(\widetilde{\theta}) = [\psi(\widetilde{\theta}^{(\lambda)}) - \psi(\widetilde{\theta})]/[\widetilde{\theta}^{(\lambda)} - \widetilde{\theta}] + o(\lambda) = -(1 + o(1))\lambda\widetilde{\theta}/[\widetilde{\theta}^{(\lambda)} - \widetilde{\theta}]$ and, hence,

$$\widetilde{\theta} - \widetilde{\theta}^{(\lambda)} = (1 + o(1))\lambda\widetilde{\theta}/\psi'(\widetilde{\theta}). \tag{4.3}$$

Since $H^{(\lambda)}(\mathbf{x}_n, \mathbf{y}_n) \geq S_{\mathbf{z}_*}^{(\lambda)}(\mathbf{x}_n, \mathbf{y}_n) = H(\mathbf{x}_n, \mathbf{y}_n) + \lambda|\mathbf{z}_*|$ for all $\lambda$, it follows by applying Theorem 2(b) on both $H(\mathbf{x}_n, \mathbf{y}_n)$ and $H^{(\lambda)}(\mathbf{x}_n, \mathbf{y}_n)$ that

$$\begin{aligned}
\limsup_{n \to \infty} |\mathbf{z}_*|/\log n &\leq [(2/\widetilde{\theta}^{(\lambda)}) - (2/\widetilde{\theta})]/\lambda \quad \text{a.s. if } \lambda > 0, \\
\liminf_{n \to \infty} |\mathbf{z}_*|/\log n &\geq [(2/\widetilde{\theta}^{(\lambda)}) - (2/\widetilde{\theta})]/\lambda \quad \text{a.s. if } \lambda < 0,
\end{aligned} \tag{4.4}$$

and Theorem 3 follows from (4.3) by letting $\lambda \to 0$ in (4.4). $\square$

## APPENDIX

**A.1. On the convexity of $\psi_\kappa$.** Let $\theta > 0$. We can express $\psi_\kappa(\theta) = \log[\sum_k a_k \times \exp(b_k\theta)]$ with $a_k \geq 0$ and $b_k$ distinct. Let $\alpha(\theta) = \sum_k a_k \exp(b_k\theta)$. Then $\psi_\kappa'(\theta) = \alpha'(\theta)/\alpha(\theta)$ and $\psi_\kappa''(\theta) = [\alpha''(\theta)/\alpha(\theta)] - [\alpha'(\theta)/\alpha(\theta)]^2$. Let $Z$ be a discrete random variable such that $P(Z = b_k) = a_k \exp(b_k\theta)/\alpha(\theta)$. Then $\psi_\kappa''(\theta) = EZ^2 - (EZ)^2 = \text{Var}(Z) \geq 0$.

**A.2. Proof of (3.11).** By the arguments just before (3.11), it suffices to show that

$$\psi(\theta) \leq \lim_{\kappa \to \infty} \xi_{2\kappa}(\theta)/2\kappa. \tag{A1}$$

Let $f_{r,s}^{(\kappa)}(\theta) = E \exp[\theta G_\kappa(\mathbf{x}_r, \mathbf{y}_s)]$ so that $h_\kappa(\theta) = \sum_{r,s \geq \kappa} f_{r,s}^{(\kappa)}(\theta)$ and let $\ell_\kappa(\theta) = \sup_{r,s \geq \kappa} f_{r,s}^{(\kappa)}(\theta)$. Let $\varepsilon > 1$ and $\theta > 0$. By (3.8), it follows that

$$\begin{aligned}
h_\kappa(\theta) &= \sum_{r,s \leq \varepsilon^\kappa + \kappa} f_{r,s}^{(\kappa)}(\theta) + 2 \sum_{r > \varepsilon^\kappa + \kappa \text{ and } s \leq \varepsilon^\kappa + \kappa} f_{r,s}^{(\kappa)}(\theta) + \sum_{r,s > \varepsilon^\kappa + \kappa} f_{r,s}^{(\kappa)}(\theta) \\
&\leq (\varepsilon^\kappa + 1)^2 \ell_\kappa(\theta) + 2(\varepsilon^\kappa + 1) \exp(\theta\kappa K_{\max}) \sum_{k > \varepsilon^\kappa} e^{-\theta g(k)} \\
&\quad + \exp(\theta\kappa K_{\max}) \bigg\{ \sum_{k > \varepsilon^\kappa} e^{-\theta g(k)} \bigg\}^2.
\end{aligned} \tag{A2}$$

Since $g(k)/\log k \to \infty$, it follows that for any $\lambda > 1$, $g(k) > (\lambda \log k)/\theta$ for all $k > \varepsilon^\kappa$ when $\kappa$ is large and

$$\text{(A3)} \quad \sum_{k > \varepsilon^\kappa} \exp[-\theta g(k)] < \sum_{k > \varepsilon^\kappa} \exp(-\lambda \log k) < \int_{\varepsilon^\kappa/2}^\infty x^{-\lambda} \, dx$$
$$= (\lambda - 1)^{-1} (\varepsilon^\kappa/2)^{-\lambda+1}.$$

Let $K_{\min} = \min_{a,b \in \mathcal{A}} K(a,b)$. Since $\ell_\kappa(\theta) \geq f_{\kappa,\kappa}^{(\kappa)}(\theta) \geq \exp(\theta \kappa K_{\min})$, it follows by choosing $\lambda > \theta(K_{\max} - K_{\min})/\log \varepsilon$ that the second and third terms on the right-hand side of (A2) are dominated by the first term as $\kappa \to \infty$. Moreover, as

$$\text{(A4)} \quad G_\kappa(\mathbf{x}_r^{(1)}, \mathbf{y}_s^{(1)}) + G_\kappa(\mathbf{x}_s^{(2)}, \mathbf{y}_r^{(2)}) \leq G_{2\kappa}(\mathbf{x}_r^{(1)} \mathbf{x}_s^{(2)}, \mathbf{y}_s^{(1)} \mathbf{y}_r^{(2)}),$$

it follows that $f_{r,s}^{(\kappa)}(\theta) f_{s,r}^{(\kappa)}(\theta) \leq f_{r+s,r+s}^{(2\kappa)}(\theta)$ for all $r, s \geq \kappa$ and, hence, by taking supremum over $r$ and $s$, we can conclude that $[\ell_\kappa(\theta)]^2 \leq \exp[\xi_{2\kappa}(\theta)]$. Hence, by (A2), (A3) and the arguments above,

$$\text{(A5)} \quad \psi(\theta) = \lim_{\kappa \to \infty} [\log h_\kappa(\theta)]/\kappa \leq 2 \log \varepsilon + \lim_{\kappa \to \infty} \xi_{2\kappa}(\theta)/2\kappa.$$

(A1) follows by letting $\varepsilon \to 1$ in (A5).

**A.3. On the existence of $\psi'(\widetilde{\theta})$.** Fix a gap penalty $g$ such that $g(k)/\log k \to \infty$ and let $\mathcal{K}$ denote the space of all symmetric matrices on $\mathcal{A} \times \mathcal{A}$ such that $K_{\max} > 0$ and $\psi(\theta) = 0$ has a unique positive solution $\widetilde{\theta}$. Induce a measure on $\mathcal{K}$ via the Lebesgue measure on the upper triangular entries of $K$. Let $\mathcal{L} = \{K \in \mathcal{K} : \psi'(\widetilde{\theta}) \text{ does not exists}\}$. We shall now show that $\mathcal{L}$ has measure zero. Consider the equivalence relation $K_1 \sim K_2$ if there exists $\lambda \in \mathbf{R}$ such that

$$\text{(A6)} \quad K_1(a,b) = K_2(a,b) + \lambda \quad \text{for all } a, b \in \mathcal{A}.$$

Let the superscript $K$ be used to signify the score matrix used. If (A6) holds, then $\psi^{(K_1)}(\theta) = \psi^{(K_2)}(\theta) + \lambda \theta$ [see line before (4.2)] and, hence, $\psi^{(K_1)}$ has a well-defined derivative at $\theta$ if and only if $\psi^{(K_2)}$ has a well-defined derivative at $\theta$. By the convexity of $\psi$, there are countably many members in each equivalence class such that $\psi'(\widetilde{\theta})$ is not well defined. If $\mathcal{L}$ is measurable, then a direct application of Fubini's theorem would show that $\mathcal{L}$ has measure 0. To show that $\mathcal{L}$ is measurable, define the distance measure $\|K - K^*\| = \max_{a,b \in \mathcal{A}} |K(a,b) - K^*(a,b)|$. Then by the convexity of $\psi^{(K)}$,

$$\text{(A7)} \quad \mathcal{L} = \bigcup_{\delta > 0} \bigcap_{\varepsilon > 0} \{K \in \mathcal{K} : \widetilde{\theta}^{(K)} > \varepsilon,$$
$$\varepsilon^{-1}[\psi^{(K)}(\widetilde{\theta}^{(K)} + \varepsilon) + \psi^{(K)}(\widetilde{\theta}^{(K)} - \varepsilon)] > \delta\},$$

where $\delta, \varepsilon$ varies over $1, \frac{1}{2}, \frac{1}{3}, \ldots$. Since $|G_\kappa^{(K)}(\mathbf{x}_m, \mathbf{y}_n) - G_\kappa^{(K^*)}(\mathbf{x}_m, \mathbf{y}_n)| \leq \kappa \|K - K^*\|$ for all $\kappa$, it follows that

(A8) $$|\psi^{(K)}(\theta) - \psi^{(K^*)}(\theta)| \leq \theta \|K - K^*\|.$$

By (A8), both $\psi^{(K)}$ and $\widetilde{\theta}^{(K)}$ are continuous with respect to $K$ and, hence, $\psi^{(K)}(\widetilde{\theta}^{(K)} + \varepsilon)$, $\psi^{(K)}(\widetilde{\theta}^{(K)} - \varepsilon)$ are also continuous with respect to $K$. The sets defined in (A7) are open and $\mathcal{L}$ is measurable.

**Acknowledgments.** The author would like to thank Amir Dembo, Wei-Liem Loh, David Siegmund, an Associate Editor and two referees for their useful references and comments.

## REFERENCES

[1] ALTSCHUL, S. F. (1998). Generalized affine gap costs for protein sequence alignment. *Proteins* **32** 88–96.
[2] ALTSCHUL, S. F., GISH, W., MILLER, W., MYERS, E. W. and LIPMAN, D. J. (1990). Basic local alignment search tool. *J. Mol. Biol.* **215** 403–410.
[3] ARRATIA, R. and WATERMAN, M. S. (1994). A phase transition for the score in matching random sequences allowing deletions. *Ann. Appl. Probab.* **4** 200–225. MR1258181
[4] BENNER, S., COHEN, M. and GONNET, G. (1993). Empirical and structural models for insertions and deletions in the divergent evolution of proteins. *J. Mol. Biol.* **229** 1065–1082.
[5] BUNDSCHUH, R. and HWA, T. (2002). Statistical mechanics of secondary structures formed by random RNA sequences. *Phys. Rev. E* **65** 031903.
[6] CHAN, H. P. (2003). Upper bounds and importance sampling of $p$-values for DNA and protein sequence alignments. *Bernoulli* **9** 183–199. MR1997026
[7] DEMBO, A., KARLIN, S. and ZEITOUNI, O. (1994). Critical phenomena for sequence matching with scoring. *Ann. Probab.* **22** 1993–2021. MR1331213
[8] DEMBO, A., KARLIN, S. and ZEITOUNI, O. (1994). Limit distribution of maximal non-aligned two-sequence segmental score. *Ann. Probab.* **22** 2022–2039. MR1331214
[9] DEWEY, T. G. (2001). A sequence alignment algorithm with an arbitrary gap penalty function. *J. Comp. Biol.* **8** 177–190.
[10] DURBIN, R., EDDY, S., KROGH, A. and MITCHISON, G. (1998). *Biological Sequence Analysis. Probabilistic Models of Proteins and Nucleic Acids.* Cambridge Univ. Press.
[11] GOTOH, O. (1982). An improved algorithm for matching biological sequences. *J. Mol. Biol.* **162** 705–708.
[12] GU, X. and LI, W. (1995). The size distribution of insertions and deletions in human and rodent pesudogenes suggests the logarithmic gap penalty for sequence alignment. *J. Mol. Evol.* **40** 464–473.
[13] GUSFIELD, D. (2001). Available at www.cs.ucdavis.edu/˜gusfield/xparall/.
[14] MILLER, W. and MYERS, E. W. (1988). Sequence comparison with concave weighting functions. *Bull. Math. Biol.* **50** 97–120. MR952421
[15] MOTT, R. (1999). Local sequence alignments with monotonic gap penalties. *Bioinf.* **15** 455–462.

[16] MOTT, R. and TRIBE, R. (1999). Approximate statistics of gapped alignments. *J. Comp. Biol.* **6** 91–112.
[17] SMITH, T. and WATERMAN, M. S. (1981). Identification of common molecular subsequences. *J. Mol. Biol.* **147** 195–197.
[18] TURNER, D. H., SUGIMOTO, N. and FREIER, S. M. (1988). RNA structure prediction. *Ann. Rev. Biophysics Biophysical Chem.* **17** 167–192.
[19] WATERMAN, M. S. (1984). Efficient sequence alignment algorithms. *J. Theor. Biol.* **108** 333–337. MR755240
[20] ZHANG, Y. (1995). A limit theorem for matching random sequences allowing deletions. *Ann. Appl. Probab.* **5** 1236–1240. MR1384373

DEPARTMENT OF STATISTICS
AND APPLIED PROBABILITY
NATIONAL UNIVERSITY OF SINGAPORE
SINGAPORE 119260
REPUBLIC OF SINGAPORE
E-MAIL: stachp@nus.edu.sg